\documentclass[11pt]{article}
\usepackage{natbib}
\usepackage{float} %Ojo, hay que cargarlo antes de hyperref que va con marsden_article, sino da problemas
\usepackage{marsden_article}
\usepackage[all]{xy}
\usepackage{amssymb}
\usepackage{enumerate,enumitem}
\usepackage{tikz-cd}
\usetikzlibrary{decorations.pathmorphing}
\usepackage[title]{appendix}
\usepackage{verbatim,xcolor,bold-extra,subcaption}
\graphicspath{ {./images/} }

\newtheoremstyle{obs}% name
  {3pt}%      Space above
  {3pt}%      Space below
  {}%         Body font
  {}%         Indent amount (empty = no indent, \parindent = para indent)
  {\bfseries}% Thm head font
  {.}%        Punctuation after thm head
  {.5em}%     Space after thm head: " " = normal interword space;
  {}%         Thm head spec (can be left empty, meaning `normal')

\theoremstyle{obs}

\newtheorem{remark}[theorem]{Remark}
\newtheorem{example}[theorem]{Example}

\newcommand{\lp}{\left(}

\newcommand\map[3]{#1\ \colon\ #2\longrightarrow#3}

\newcommand{\st}{\;\ifnum\currentgrouptype=16 \middle\fi|\;}

\usepackage{comment}
 \newcommand{\lvec}[1]{\overleftarrow{#1}}
\newcommand{\rvec}[1]{\overrightarrow{#1}}

\DeclareMathOperator{\Tr}{Tr}

\setenumerate[1]{label=(\roman*), ref=(\roman*)}
\usepackage[all]{xy}% Para los diagramas.

\newtheoremstyle{obs}% name
  {3pt}%      Space above
  {3pt}%      Space below
  {}%         Body font
  {}%         Indent amount (empty = no indent, \parindent = para indent)
  {\bfseries}% Thm head font
  {.}%        Punctuation after thm head
  {.5em}%     Space after thm head: " " = normal interword space;
  {}%         Thm head spec (can be left empty, meaning `normal')
\theoremstyle{obs}

\newcommand{\R}{\mathbb{R}}      %Numeros reales
\begin{document}
\title{Lie-Poisson integrators}
\author{
\textsc{David Mart\'{\i}n de Diego}\thanks{david.martin@icmat.es} \\
\small
Instituto de Ciencias Matem\'aticas (CSIC-UAM-UC3M-UCM) \\
 \small  C/Nicol\'as
Cabrera 13-15, 28049 Madrid, Spain
}

\maketitle

\begin{abstract}

In this paper, we discuss the geometric integration of hamiltonian systems on Poisson manifolds, in particular, in the case, when the Poisson structure is induced by a Lie algebra, that is, it is a Lie-Poisson structure.

A Hamiltonian system on a Poisson manifold $(P, \Pi)$ is a smooth manifold $P$ equipped with a bivector field $\Pi$ satisfying $[\Pi, \Pi\;]=0$ (Jacobi identity), inducing the Poisson   bracket on $C^{\infty}(P)$, $\{f, g\}\equiv \Pi(df, dg)$ where $f, g\in C^{\infty}(P)$. For any $f\in C^{\infty}(P)$ the Hamiltonian vector field is defined by $X_f(g)=\{g, f\}$. The Hamiltonian vector fields $X_f$ generate an integrable generalized distribution on $P$ and the leaves of this foliation are symplectic. The flow of any hamiltonian vector field preserves the Poisson structure, it fixes each leaf and the hamiltonian itself is a first integral. 

It is important to characterize numerical methods preserving some of these fundamental properties of the hamiltonian flow on Poisson manifolds (geometric integrators). We discuss the difficulties of deriving these Poisson methods using standard techniques and we present some modern approaches to the problem.

\vspace{3mm}


\vspace{3mm}

\end{abstract}

\tableofcontents

\section{Introduction}

Many applied  dynamical systems often display a variety of geometric structures in their mathematical description. For instance, you can think in the standard autonomous mechanical systems given by  kinetic energy minus potential energy (that is, a mechanical Lagrangian system) (see \cite{foundation} and references therein).  It  is easy to show that  these classes of systems are preserving at least a  constant of the motion, the total energy, and moreover they also preserve a volume form or, even more precisely, a  symplectic structure. Other constants of motion may arise if we are in presence of  additional symmetries of the Lagrangian as, for example, the preservation of the angular momentum in the case of mechanical systems in a central force field (see \cite{marsden3,MR818988}).

As a short introduction, the Lagrangian formulation of mechanics can be based on the variational principles founded  by Euler and Lagrange given a mathematical interpretation 
of Newton's fundamental law of force balance $F = ma$. 
The first element is the 
configuration space $Q$, which describes all the possible positions of the mechanical system. In more mathematical terms, the configuration space is geometrically described   by a  manifold, finite or infinite dimensional. For clarity on the exposition, we will assume in the sequel that our configuration space is  finite dimensional, $\dim Q=n$, and locally described by a set of coordinates  denoted by $(q^i)$, $i = 1, \ldots,n$, that prescribe the admissible configurations of the system under study. 
Since the Newton's fundamental law implies a system of second  order differential equations, then the dynamics of a mechanical system is described by a special vector field on the phase space of positions and velocities. Geometrically, one must introduce the  
 velocity phase space $TQ$, the tangent bundle of $Q$, with induced coordinates 
$(q^i, \dot{q}^i)$, $i = 1, \ldots, n$ and canonical projection $\tau_Q: TQ\rightarrow Q$, given by $\tau_Q(q^i, \dot{q}^i)=(q^i)$. Therefore, any vector $v_q\in T_qQ$
is expressed in local coordinates as
\[
v_q=\left.v^i\frac{\partial}{\partial q^i}\right |_{q}
\]
and  an arbitrary  second order differential equation  is now geometrically described as an special class of vector field of $TQ$:
\[
\Gamma=\dot{q}^i\frac{\partial}{\partial q^i}+\Gamma^i(q, \dot{q})\frac{\partial}{\partial \dot{q}^i}\; .
\]

One of the most important properties of variational calculus is that the dynamics is completely specified by an unique  function, the lagrangian $L: TQ\rightarrow {\mathbb R}$. In mechanical systems, $L$ is
the kinetic minus the potential energy of the system.
 Hamilton's principle states that the variation
of the action is stationary:
\[
\delta \int_{0}^h L(q(t), \dot{q}(t))\; dt=0\; .
\]
In this principle, one chooses curves $q: [0, h]\rightarrow Q$ 
joining two fixed points in $Q$, $q(0)=q_0$ and $q(h))=q_1$. The extremals are characterized  by 
the solutions of the Euler-Lagrange equations:
\[
\frac{d}{dt}\left(\frac{\partial L}{\partial \dot{q}^i}\right)-\frac{\partial L}{\partial q^i}=0, \  1\leq i\leq n
\]
which is a  implicit system of second order differential equations. 
For Lagrangians that are purely kinetic energy, then the Lagrangian is completely defined by a riemannian metric and the solutions of the Euler-Lagrange
equations are the riemannian  geodesics. 

Of course these equations can be generalized in many ways, for instance, it is possible to add
external forces introduced  on the right
hand side of the  Euler-Lagrange equations (see, for instance, \cite{BKMR96}); but also is possible to add constraints introducing the study of the so-called nonholonomic and vakonomic mechanics (see \cite{arnoldrkfufaev,arnoldsist3}). 

In many cases is interesting to  pass to the Hamiltonian formalism introducing the conjugate momenta
$p_i=\frac{\partial L}{\partial \dot{q}^i}$, $1\leq i\leq n$. 
Geometrically, it corresponds to change the tangent bundle $TQ$ by its dual bundle $T^*Q$ with coordinates $(q^i, p_i)$, $1\leq i\leq n$, given by the Legendre transformation
which in canonical coordinates is defined by 
\[
Leg_L(q^i, p_i)=(q^i, \frac{\partial L}{\partial \dot{q}^i})\; .
\]
The Lagrangian is called regular when this change of variables is invertible, then  we can express $\dot{q}^i=\dot{q}^i(q, p)$ and introduce the Hamiltonian
\[
H(q^i, p_i)=p_j \dot{q}^j(q, p) - L(q^i, \dot{q}^i(q, p))\; .
\]
One shows that the Euler-Lagrange equations are equivalent to Hamilton's equations:
\begin{eqnarray*}
\frac{ d q^i}{dt}&=&\frac{\partial H}{\partial p_i}\\
 \frac{ d p_i}{dt}&=&-\frac{\partial H}{\partial q^i}
\end{eqnarray*}
with $i=1, \ldots , n$. 
In $T^*Q$ it is more transparent to explore the qualitative properties of mechanical systems. For instance, on $T^*Q$ is defined the canonical symplectic 2-form
\[
\omega_Q=dq^i\wedge dp_i\, ,
\]
and the Hamiltonian equations are written as 
\[
i_{X_H}\omega_Q=dH\; .
\]
The integral curves of $X_H$ are solutions of the Hamilton's equations.
This implies that the flow $\Phi: D\subseteq {\mathbb R}\times T^*Q\rightarrow T^*Q$  of $X_H$ verifies $\Phi_t^*\omega_Q=\omega_Q$, that is, the flow is a symplectomorphism for all $t$.

These Hamilton's equations can be rewritten in terms of the canonical Poisson bracket form as
\[
\frac{df}{dt}=\{H, f\}
\]
where the  Poisson bracket is locally expressed as
\[
\{f,g\}=\sum_{i=1}^n \left(
\frac{\partial f}{\partial p_i}
\frac{\partial g}{\partial q^i}
-
\frac{\partial f}{\partial q^i}
\frac{\partial g}{\partial p_i}\right)\; ,\quad f, g\in C^{\infty}(T^*Q)
 \]

Of course, since typically is impossible to explicitly integrate a given dynamical systems it is useful to use computational schemes preserving as much as possible the geometric structures inherently associated to the original system, mainly if you are interested in long-term predictions. See \cite{hairer,MR1270017,blanes}.

A good example of these geometric integrators are symplectic schemes for a hamiltonian system determined by $H: \R^{2n}\rightarrow \R$. The Hamilton equations of motion are  written as
\begin{equation}\label{asd}
\left(
\begin{array}{c}
\dot{q}\\
\dot{p}
\end{array}
\right)
=
\left(
\begin{array}{cc}
0&{\mathbb I}\\
-\mathbb{I}&0
\end{array}
\right)
\left(\begin{array}{c}
\frac{\partial H}{\partial q}(q,p)\\
\frac{\partial H}{\partial p}(q, p)
\end{array}
\right)
\end{equation}
A second-order symplectic integrator for this system is given by the mid-point rule
\begin{equation}\label{asd}
\left(
\begin{array}{c}\displaystyle \frac{q_{k+1}-q_k}{h}\\
\displaystyle\frac{p_{k+1}-p_k}{h}
\end{array}
\right)
=
\left(
\begin{array}{cc}
0&{\mathbb I}\\
-\mathbb{I}&0
\end{array}
\right)
\left(
\begin{array}{c}
\displaystyle\frac{\partial H}{\partial q}\tiny\left(\frac{q_{k}+q_{k+1}}{2},\frac{p_{k+1}+p_{k}}{2}\right) \\
\displaystyle\frac{\partial H}{\partial p}\left(\frac{q_{k}+q_{k+1}}{2},\frac{p_{k+1}+p_{k}}{2}\right)
\end{array}
\right)
\end{equation}
The mid-point rule was studied by \cite{Krishnaprasad} for the more general situation of hamiltonian systems defined by Poisson manifolds (see Section \ref{Poisson}) where now the structure is defined by a Hamiltonian function $H: \R^m\rightarrow \R$ and the Poisson tensor is determined by a matrix 
$\Pi=(\Pi^{ij})_{1\leq i, j\leq m}$ satisfying
\begin{itemize}
\item Skew-symmetry $\Pi=-\Pi^T$.
\item 
Jacobi identity 
\[
0=\sum_{l=1}^m\left(
\frac{\partial \Pi^{ij}}{\partial z^l}\Pi^{lk}+\frac{\partial \Pi^{jk}}{\partial z^l}\Pi^{li}+\frac{\partial \Pi^{ki}}{\partial z^l}\Pi^{lj}\right), \quad i, j, k=1, \ldots, m
\]
\end{itemize}
The hamiltonian equations in this more general case are written as
\begin{equation}\label{aqw}
\dot{z}=\Pi(z)\nabla H(z)
\end{equation}
A transformation $\varphi: \R^m\rightarrow \R^m$ is said to be Poisson map for  $(\R^m, \Pi)$ if 
\[
D\varphi(z)\Pi(z)D\varphi(z)^T=\Pi(\varphi(z))\; .
\]
The flow $\Phi_t$ of the dynamical system (\ref{aqw}) is a Poisson transformation and moreover due to the skew-symmetry of the matrix $\Pi$, the energy is preserved, that is
$H(\Phi_t(z))=H(z)$. Moreover, if $C: {\mathbb R}^m\rightarrow {\mathbb R}$ is a function such that $\Pi(z)\nabla C (z)=0$  (a Casimir function), then 
\[
C(\Phi_t(z))=C(z)\; .
\]
All them are important qualitative properties of the Hamiltonian vector field on a Poisson manifold.  

As an interesting example  showing the difficulties to preserve structures in cases different to the explored in \cite{Krishnaprasad},  they studied the mid-point rule applied to (\ref{aqw}), that is, 
\begin{equation}\label{poi}
\frac{z_{n+1}-z_n}{h}=\Pi \left(\frac{z_{n+1}+z_n}{2}\right)\nabla H \left(\frac{z_{n+1}+z_n}{2}\right)
\end{equation}
It is trivial to check  that if $\Pi$ is constant, that is $\Pi(z)=\Pi$ constant, then the mid-point rule is a Poisson automorphism. However, in 1993, Austin, Krisnaprasad and Wang showed that the mid-point rule is an ``almost-Poisson integrator", in the sense that preserves the Poisson structure up to second-order
\[
D\Phi_H^h(z)\Pi(z)D\Phi_H^h(z)^T-\Pi(\Phi_H^h(z))=O(h^3)\; .
\]
where $\Phi_H^h: \R^n\rightarrow \R^n$ is the map implicitly defined by Equation (\ref{poi}), that is, $z_{n+1}=\Phi_H^h(z_n)$. 

In this paper, we will study some recent developments on geometric integration on Poisson manifolds and more specially in Lie-Poisson systems, that is, when the Poisson bracket is linear on a vector space. 
This particular case is specially interesting since is close related with invariant mechanical  systems defined on Lie groups, as for instance,  equations of the  rigid body, heavy top and fluids as special cases. 
His  background led to  Poincar\'e (see \cite{poincare}) and the equations are determined  giving a Lagrangian on a Lie algebra (see \cite{HMR98} for more details). The corresponding hamiltonian description is determined using the Lie-Poisson bracket as we will see.

The paper is structured as follows.  In Section 2 we will introduce some basic notions of Poisson manifolds with an spatial emphasis in Lie-Poisson brackets. In Section 3 we  will describe the Euler-Poincar\'e equations and Lie-Poisson equations. Finally, in Section 4 we will summarize  some modern methods  to numerically integrate the Lie-Poisson equations using geometric integrators, presenting some original results to this topic.

\section{Poisson manifolds}\label{Poisson}

Poisson manifolds appears as a natural generalization of symplectic manifolds. As we will see along this paper, 
Poisson manifolds occur as phase spaces for classical mechanics but moreover it is a concept relevant in quantum mechanics (see, for more information, \cite{MR723816,MR0501133,MR882548,MR960879,MR1269545,MR1747916}).

We recall that a {\bf Poisson structure} on a  differentiable manifold $P$ is given by a bilinear map
\[
\begin{array}{rcc}
C^{\infty}(P)\times C^{\infty}(P)&\longrightarrow& C^{\infty}(P)\\
(f, g)&\longmapsto& \{f, g\}
\end{array}
\]
called the {\bf Poisson bracket}, satisfying the following properties: 
\begin{itemize}
\item[(i)] \emph{Skew-symmetry},  $\{g, f\}=-\{f, g\}$;
\item[(ii)] \emph{Leibniz rule}, $\{fg, h\}=f\{g, h\}+g\{f, h\}$; 
\item[(iii)] \emph{Jacobi identity},  $\{\{f, g\}, h\}+\{\{h, f\}, g\}+\{\{g, h\}, f\}=0$;
\end{itemize}
for all $f, g, h\in C^{\infty}(P)$. 

The situation for Poisson manifolds is in some sense similar to the case of a riemannian manifold where, after one has fixed a riemannian metric, each  function determines the corresponding gradient vector field.

Given a Poisson manifold with bracket $\{\; ,\;\}$ and a function $f\in C^{\infty}(P)$ we may associate a unique vector field $X_f\in {\mathfrak X}(P)$, the {\bf Hamiltonian vector field}: 
\[
X_f(g)=\{f, g\}\; .
\]
%The assignment $f\mapstro X_f$  is a Lie algebra morphism between $C^{\infty}(P)$ and ${\mathfrak X}(P)$.

Moreover, on a Poisson manifold there exists a unique bivector field $\Pi$, a Poisson bivector  (that is, a twice contravariant skew symmetric differentiable tensor) 
such that 
\begin{equation*}
  \{f,g\}:=\Pi(df ,d g), \qquad f, g\in C^\infty (P)
\end{equation*}
The bivector field $\Pi$ is called the {\bf Poisson tensor} and the Poisson structure is usually denoted by $(P, \Pi)$. 
The Jacobi identity in terms of the bivector $\Pi$ is written as
\[
[\Pi, \Pi\,]=0
\]
where here $[\; , \; ]$ denotes the Schouten bracket.  

Typical examples of Poisson manifolds are:
\begin{itemize}
\item {\bf Symplectic structures}, that is, if $P$ is equipped with a non-degenerate
closed $2$-form
$\omega$ on $P$, then its inverse is a Poisson bivector.
\item   The {\bf dual of
a Lie algebra} ${\mathfrak g}^*$. If ${\mathfrak g}$ is a Lie algebra with Lie bracket $[\; ,\; ]$, then it is defined a   Poisson bracket on
 ${\mathfrak g}^*$  by 
 \[
 \{\xi, \eta \}(\alpha) = -\langle\alpha ,[\xi,\eta]\rangle\; ,
 \] where $\xi$ and $\eta \in
{\mathfrak g}$ are equivalently considered as linear forms on ${\mathfrak g}^*$, and $\alpha \in {\mathfrak g}^*$.
This  linear Poisson structure on ${\mathfrak g}^*$ is called the
{\it Kirillov-Kostant-Souriau Poisson structure}.
\end{itemize}

Taking coordinates  $(x^i)$, $1\leq i\leq \dim P=m$ and $\Pi^{ij}$ the components of the Poisson bivector, that is
\[
\Pi^{ij}=\{x^i, x^j\}
\]
 then if $f, g\in C^{\infty}(P)$ 
 \[
 \{f, g\}=\sum_{i,j=1}^m\{x^i, x^j\}\frac{\partial f}{\partial x^i}\frac{\partial g}{\partial x^j}= \sum_{i,j=1}^m\Pi^{ij}\frac{\partial f}{\partial x^i}\frac{\partial g}{\partial x^j}\; .
\] 
 Observe that the $m\times m$ matrix $(\Pi^{ij})$  verifies the following properties:
 \begin{itemize}
\item[(i)] \emph{Skew-symmetry}, $\Pi^{ij}=-\Pi^{ji}$  
\item[(ii)] \emph{Jacobi identity} 
\[
\sum_{l=1}^m\left(
\Pi^{il}\frac{\partial \Pi^{jk}}{\partial x^l}+\Pi^{kl}\frac{\partial \Pi^{ij}}{\partial x^l}+\Pi^{jl}\frac{\partial \Pi^{ki}}{\partial x^l}\right)=0\, ,\quad  i, j, k=1,\ldots, m.
\]
 \end{itemize}
 
Define ${\sharp}^{\Pi} : T^{*}P \to TP$ by
\[
{\sharp}^\Pi (\alpha) = \iota_\alpha \Pi=\Pi(\alpha, \cdot)  ,
\]
where $\alpha \in T^*P$,
and $\langle \beta, \iota_\alpha \Pi\rangle=\Pi(\alpha, \beta)$ for all $\beta \in T^*P$.
The rank of $\Pi$ at $p\in P$ is $\hbox{rank }{\sharp}^\Pi_p: T^*_pP\rightarrow T_pP$.  Because of the skew-symmetry of $\Pi$, we know that the rank of $\Pi$  at a point $p\in P$ is an even integer. 

Fixed a function  $H \in
C^\infty(P)$, the hamiltonian function, we have the corresponding Hamiltonian vector field:
\[
X_H = {\sharp}^\Pi (dH) .
\]
Therefore, on a Poisson manifold, a function determines the following  dynamical system:
\begin{equation}\label{hamil-eq}
\frac{dx}{dt}=X_H(x(t))\; .
\end{equation}
Moreover, a function $f\in C^{\infty}(P)$  is a first integral of the Hamiltonian vector field  $X_H$ if for any solution 
 $x(t)$ of Equation (\ref{hamil-eq}) we have
\[
\frac{df}{dt}(x(t))=0\; .
\]
In other words, if  $X_H(f) = 0$ or, equivalently, $\{H, f\} = 0$. In particular,  the
hamiltonian is a conserved quantity since $\{H, H\}=0$ by the skew symmetry of the bracket.
For any Poisson manifold $(P,\Pi)$  a function $C\in C^{\infty}(P)$
 is called a
{\bf Casimir function} of $\Pi$ if $X_C = 0$, i.e, if $\{C, g\} = 0$, for all $g\in C^{\infty}(P)$.

Additionally, it is easy to show that if $f_1$ and $f_2$ are first integrals of $X_H$ then so is $\{f_1, f_2\}$. The proof is a direct consequence of the Jacobi identity since: 
\[
\{H, \{f_1, f_2\}\}=\{\{H, f_1\}, f_2\}+\{\{f_2, H\}, f_1\}=0
\]

\begin{example}{\rm 
Consider the Lie algebra ${\mathfrak g}={\mathfrak so}(3)$ of $3\times 3$-skew symmetric matrices, that it is also identified with ${\mathbb R}^3$ with the vector product $\times$.  Also ${\mathfrak so}(3)^*$ is identified with ${\mathbb R}^3$ and the corresponding Poisson bracket of two functions $f, g\in C^{\infty}({\mathbb R}^3)$: 
\[
\{f, g\}(x,y, z)=
\left|
\begin{array}{ccc}
\frac{\partial f}{\partial x}&\frac{\partial f}{\partial y}&\frac{\partial f}{\partial z}\\
\frac{\partial g}{\partial x}&\frac{\partial g}{\partial y}&\frac{\partial g}{\partial z}\\
x&y&z
\end{array}
\right|\; .
\]
Observe that the function $C=\frac{1}{2}(x^2+y^2+z^2)$ is  a Casimir function for this bracket.

Now,  considering the typical hamiltonian for the rigid body   $$H(x,y,z)=\frac{1}{2}\left(\frac{x^2}{I_1}+\frac{y^2}{I_2}+\frac{z^2}{I_3}\right),$$ then the hamiltonian vector field $X_h$ gives the following system of equations:  
\begin{eqnarray*}
\dot{x}&=&\{H, x\}=\frac{I_2-I_3}{I_2I_3}yz\\
\dot{y}&=&\{H, y\}=\frac{I_1-I_3}{I_2I_3}xz\\
\dot{z}&=&\{H,z\}=\frac{I_1-I_2}{I_1I_2}xy \hspace{5cm}\diamond
\end{eqnarray*}
}
\end{example}

One can show that $\hbox{Im} (\sharp^{\Pi})\subseteq TP$ is an involutive generalized distribution called the {\bf characteristic distribution}. Then,  given a Poisson structure $\Pi$ we have that the differentiable manifold  $P$ is foliated by leaves (that is, immersed
submanifolds of varying dimensions) such that their  tangent spaces are  given by
$\hbox{Im } (\sharp_p^\Pi)$.
Fixed a point $p\in P$ the kernel of 
${\sharp}^{\Pi}_p: T^*_pP\rightarrow \hbox{Im } \sharp^\Pi_p$ is precisely the annihilator 
\[
(\hbox{Im } \sharp_p^\Pi)^0=\{\alpha_p\in T^*_pP\; |\; 
\langle \alpha_p,  X_p\rangle=0\, ,  \forall X_p\in \hbox{Im }{\sharp}^\Pi_p\}\; ,
\]
 We derive an isomorphism 
\[
T^*_pP/ (\hbox{Im } \sharp_p^\Pi)^0\equiv \hbox{Im } (\sharp^\Pi_p)^*\longrightarrow \hbox{Im } (\sharp^\Pi_p)\; ,
\]
and, in consequence,  a linear symplectic form $\omega_p$ on $\hbox{Im } (\sharp^\Pi_p)$. 
This shows that each leaf is equipped with a symplectic form, that is, a Poisson structure defines a symplectic foliation on $P$.  
 
 \begin{theorem}{\bf Darboux-Weinstein coordinates.}
 Let $p\in P$ be an arbitrary point in a Poisson manifold $(P, \Pi)$ of $\dim P=m$ with rank of $\Pi$ at $p$ is $2k$. 
 There is a system of local coordinates
$(U, (q^1, \ldots, q^k, p_1, \ldots, p_k, y^1, \ldots, y^s))$ centered at $p$  such that: 
\[
\{f, g\}=
\sum_{i=1}^k\left(
\frac{\partial f}{\partial p_i}\frac{\partial g}{\partial q^i}-\frac{\partial f}{\partial q^i}\frac{\partial g}{\partial p_i}\right)
+\sum_{i<j}^s \phi^{ij}(y)\frac{\partial f}{\partial y^i}\frac{\partial g}{\partial y^j}
\]
where $\phi^{ij}=\phi^{ij}(y)$ are functions that depend only on the $(y^1, \ldots, y^s)$
and vanish at $p$.
 \end{theorem}
Observe that locally
a Poisson bracket splits into two pieces: an standard Poisson bracket
on $\R^{2k}$
and a singular Poisson bracket vanishing at $p$.

In particular, if $(P,\Pi)$ is a Poisson manifold of constant rank $2k$, then for every point $p\in M$, there exists a system $(U, (q^1, \ldots, q^k, p_1, \ldots, p_k, y^1, \ldots, y^s))$ of local coordinates  
 such that: 
\[
\{f, g\}=
\sum_{i=1}^k\left(
\frac{\partial f}{\partial p_i}\frac{\partial g}{\partial q^i}-\frac{\partial f}{\partial q^i}\frac{\partial g}{\partial p_i}\right)\; .
\]

Let $(P_1, \Pi_1)$ and $(P_2, \Pi_2)$ be two Poisson manifolds and $\varphi: P_1\rightarrow P_2$ be a differentiable map. We say that $\varphi$ is a {\bf Poisson morphism} if for all $f, g\in C^{\infty}(P_2)$, 
\[
\varphi^*\{f, g\}=\{\varphi^* f, \varphi^*g\}\; .
\]
It is easy to show that this property is equivalent to
\[
T_x\varphi( X_{\varphi^*f}(x))=X_f(\varphi(x))\; ,
\]
for all $f\in C^{\infty}(P_2)$ and $x\in P_1$.  
This immediately implies that if $\sigma: I\subseteq \R\rightarrow P_1$ is an integral curve of $X_{\varphi^*f}$ then 
$\varphi\circ \sigma: I\subseteq \R\rightarrow P_2$ is an integral curve of $X_f$

Let $(P, \Pi)$ be a Poisson manifold. A submanifold $i_N: N\hookrightarrow P$ is a {\bf Poisson submanifold} if we can define a Poisson bracket $\{\; ,\; \}_N$ by 
\[
\{\tilde{f}, \tilde{g}\}_N=\{f, g\}_{|N}\, ,
\]
where $f, g\in C^{\infty}(P)$ are arbitrary smooth extensions of $\tilde{f}, \tilde{g}\in  C^{\infty}(N)$. Then, the inclusion $i_N$ is a Poisson morphism. It is possible to show that $N$ is a Poisson submanifold if and only if for each each point $p\in N$ 
\[
\hbox{Im } (\sharp_p^\Pi)\subseteq T_pN\; .
\]
In other words, if all the Hamiltonian vector fields $X_H\in {\mathfrak X}(P)$ are tangent to  $N$.

A submanifold $N$ of a Poisson manifold $(M, \Pi)$  is called {\bf coisotropic}
if
\[
 (\sharp^\Pi)(TN^0)\subseteq TN
\]
or equivalently if the ideal 
$I_N =\left\{f\in C^{\infty}(P) \; | \; f_{|N} = 0\right\}$  is closed under the Poisson bracket $\{\;,\; \}$, that is, 
\[
\{I_N, I_N\}\in I_N\; .
\]
Two interesting examples are precisely the extreme cases: 
\begin{itemize}
\item $ (\sharp^\Pi)(TN^0)=0$, these are exactly the {\bf Poisson submanifolds}; 
\item  $(\sharp^\Pi)(TN^0)=TN$, these submanifolds are called {\bf Lagrangian submanifolds}.
\end{itemize}

Let $(P_1, \Pi_1)$ and $(P_2, \Pi_2)$ be Poisson manifolds, then
\begin{itemize}
\item If $N\subset P_1$ is a coisotropic submanifold, $f: P_1\rightarrow P_2$ a Poisson map and $f(N)$ is a submanifold of $P_2$ , then $f(N)$ is a 
coisotropic submanifold of $(P_2, \Pi_2)$;
\item Denote by $(\overline{P_1}\times {P_2}, \Pi_{\overline{P_1}\times {P_2}})$, the Poisson manifold   on $P_1\times P_2$ where the Poisson structure $\Pi_{\overline{P_1}\times {P_2}}$ is given by:
\[
\Pi_{\overline{P_1}\times {P_2}}(\alpha_1+\alpha_2, \beta_1+\beta_2)=\Pi_2(\alpha_2, \beta_2)-\Pi_1(\alpha_1, \beta_1) \; ,
\]
for any $\alpha_1, \beta_1\in T^*P_1$ and  $\alpha_2, \beta_2\in T^*P_2$. Then, it is easy to check that $\varphi: P_1\rightarrow P_2$ is a Poisson map if and only if 
\[
\hbox{Graph}(\varphi)=\{(p_1, \varphi(p_1)\; |\, p_1\in P_1\}
\]
is a coisotropic submanifold of  $(\overline{P_1}\times {P_2}, \Pi_{\overline{P_1}\times {P_2}})$.
\end{itemize}
This last property is important for constructing geometric integrators for Poisson systems.

A {\bf symplectic realization} of a Poisson manifold $(P, \Pi)$ is a Poisson map $J: M\rightarrow P$ from a symplectic manifold $(M, \omega)$ to $P$.
 If $J$ is a surjective submersion we say that it is a {\bf full symplectic realization}.
We have that $J: M\rightarrow P$ maps local Lagrangian submanifolds to coisotropic submanifols of $(P, \pi)$.  In particular recovering the graphs of Poisson morphisms. 
If $(M, \omega) = (T^*{\mathbb R}^n, dq^i\wedge dp_i)$ then  $(q^i, p_i)$ are called canonical or Clebsch variables for $(P, \Pi)$. 
See \cite{MR3215845} for the application of this concept to design Poisson integrators from symplectic integrators.

Another interesting geometric concept is the notion of {\bf dual pair}.  A dual pair is a pair of smooth Poisson maps $(\Psi_1: M\rightarrow P_1$, $\Psi_2: M\rightarrow P_2)$ where  $(M, \omega)$ is a symplectic manifold  and  $(P_1, \Pi_1)$ and $(P_2, \Pi_2)$ are Poisson manifolds, such that for each point  $x\in M$, then: 
\[
\ker(T\Psi_1) =
\ker(T\Psi_2)^{\perp, \omega}
\qquad
\ker(T\Psi_2) = 
\ker(T\Psi_1)^{\perp, \omega}
\]
where
\[
 \ker(T\Psi_i)^{\perp, \omega}=
 \{ V\in TM	\; |\; \omega(V, 	\ker(T\Psi_i))=0\}, \quad i=1,2\; .
 \]
In many cases, the Poisson maps $\Psi_i$ are momentum mappings associated to Lie algebra actions on $M$. 
For this reason dual pairs  are a fundamental ingredient to construct Poisson integrators (see \cite{MR1115869,MR967725,MR3686003}).

\begin{example}{\bf Action of a Lie group on its cotangent bundle.}{
\rm  
Consider a Lie group $G$ and its cotangent bundle  $\pi_G: T^*G\rightarrow G$. For  each $g$ and $h$ in $G$ we have the corresponding  left and right translations 
\[
{\mathcal L}_g(h)=gh\; ,\quad  {\mathcal R}_g(h)=hg
\]

Denote by $\hat{{\mathcal L}}_g: T^*G\rightarrow T^*G$ and $\hat{{\mathcal R}}_g: T^*G\rightarrow T^*G$ the canonical lifts of ${\mathcal L}_g$ and ${\mathcal R}_g$, respectively,  defined by
\begin{eqnarray*}
\langle \hat{{\mathcal L}}_g(\mu_h), X_{gh}\rangle&=&\langle \mu_h, T{\mathcal L}_{g^{-1}} X_{gh}\rangle\\
 \langle \hat{{\mathcal R}}_g(\mu_h), Y_{hg}\rangle&=&\langle \mu_h, T{\mathcal R}_{g^{-1}} Y_{hg}\rangle
\end{eqnarray*}
 These two actions are Hamiltonian and have as momentum maps,
respectively, the maps $J_L: T^*G\rightarrow {\mathfrak g}^*$ and $J_R: T^*G\rightarrow {\mathfrak g}^*$ defined by 
\begin{eqnarray*}
J_L(\mu_g)&=&\hat{{\mathcal R}}_{g^{-1}}(\mu_g)\\
J_R(\mu_g)&=&\hat{{\mathcal L}}_{g^{-1}}(\mu_g)
\end{eqnarray*}
where ${\mathfrak g}$ is the Lie algebra of $G$.

As we have said on ${\mathfrak g}^\ast$ we have  the bracket: 
\begin{eqnarray*}
\{f, g\}(\mu)&=&-
\langle \mu, [df(\mu), dg(\mu)]\rangle
=-\langle \mu, \left[ \frac{\delta f}{\delta \mu}, \frac{\delta g}{\delta \mu}\right]\rangle
\end{eqnarray*}
where $\mu\in {\mathfrak g}^\ast$ and $\frac{\delta f}{\delta \mu}: {\mathfrak g}^\ast\rightarrow {\mathfrak g}$ stands for the functional derivative of $f$ and where $[\, ,\,  ]$ is the Lie algebra bracket on ${\mathfrak g}$. 

The map $J_L: T^*G\rightarrow {\mathfrak g}^*$ is a Poisson map where $T^*G$ is equipped 
with the canonical Poisson bracket and ${\mathfrak g}^*$ with the Poisson structure induced by 
$\{\; , \;\}$. Also,  $J_R: T^*G\rightarrow {\mathfrak g}^*$ is a Poisson map between $(T^*G, \omega_G)$ and 
$({\mathfrak g}^*, -\{\; , \;\})$.

The pair of maps $(J_L, J_R)$ is are both surjective submersions. It is the simpler example of a dual pair and both are examples of  symplectic realizations of a Poisson manifold.

}\end{example}

\section{Euler-Lagrange equations when the configuration  is a Lie group}

Let $G$ be a Lie group. 
The left multiplication ${\mathcal L}_g$ allows us to trivialize the tangent bundle $TG$ and the cotangent bundle $T^*G$ as follows
\begin{eqnarray*}
TG&\to&G\times {\mathfrak g}\, ,\qquad (g, \dot{g})\longmapsto (g, g^{-1}\dot{g})=(g, T_g{\mathcal L}_{g^{-1}}\dot g)=(g, \xi)\; ,\\
T^*G&\to&G\times {\mathfrak g}^*,\qquad  (g, \alpha_g)\longmapsto (g, T^*_e {\mathcal L}_g(\alpha_g))=(g, \alpha)\; ,
\end{eqnarray*}
 where ${\mathfrak g}=T_eG$ is the Lie algebra of $G$ and $e$ is the neutral element of $G$.

Given a lagrangian $L: TG\rightarrow {\mathbb R}$ then in left trivialized coordinates the classical Euler-Lagrange equations are given by
\begin{eqnarray}
\frac{d}{dt}\left(\frac{ \delta L}{\delta \xi}\right)&=&ad_{\xi}^*\frac{\delta L}{\delta\xi}+T^*_e{\mathcal L}_g\left(\frac{\delta L}{\delta\xi}\right)\label{elt}\\
\dot{g}&=&g\xi
\end{eqnarray}

Therefore, if $L: TG\rightarrow {\mathbb R}$ is left invariant we can define the reduced lagrangian $l: {\mathfrak g}\rightarrow {\mathbb R}$ by
\[
 l(\xi)=L(e, \xi)
 \]
 that is, $l$ is the restriction of $L$ to ${\mathfrak g}$. In this case, the corresponding Euler-Lagrange equations 
 \[
 \frac{d}{dt}\left(\frac{ \delta l}{\delta \xi}\right)=ad_{\xi}^*\frac{\delta l}{\delta\xi}
\]
 are known as the (left-invariant) {\bf Euler-Poincar\'e equations} (see \cite{marsden3,MR2847777,MR2548736}). 
From these equations we can reconstruct a solution $t\mapsto g(t)$ of the Euler-Lagrange equations with initial condition $g(0)=g_0$ and $\dot{g}(0)=v_0$ as follows. First we solve the first order differential equation 
  \[
 \frac{d}{dt}\left(\frac{ \delta l}{\delta \xi}\right)=ad_{\xi}^*\frac{\delta l}{\delta\xi}
\]
with initial condition $\xi(0)=g_0^{-1}v_0$ and then with the solution $t\rightarrow \xi(t)$ we must   solve the reconstruction equation: 
\[
\dot{g}(t)=g(t)\xi(t) \hbox{   with   } g(0)=g_0\; .
\]

\subsection{Euler-Arnold equations} (See, for more details, \cite{MR3242761})

Analogously, using this left trivialization it is possible to write the classical Hamiltonian equations for a Hamiltonian function $H: T^*G\rightarrow \R$
using  a different and interesting perspective.

For instance, it is easy to show that the canonical structures of the cotangent bundle: the Liouville 1-form $\theta_G$ and the canonical symplectic 2-form $\omega_G$, are now rewritten using this left trivialization as follows:
\begin{eqnarray}
(\theta_G)_{(g, \alpha)}(\xi_1, \nu_1)&=&\langle \alpha, \xi_1\rangle\; ,\\
(\omega_G)_{(g, \alpha)}\left( (\xi_1, \nu_1), (\xi_2, \nu_2)\right)&=&-\langle \nu_1, \xi_2\rangle + \langle \nu_2, \xi_1\rangle+\langle\alpha, [\xi_1, \xi_2]\rangle\; ,\label{omega}
\end{eqnarray}
with $(g, \alpha)\in G\times {\mathfrak g}^*$,
where $\xi_i\in {\mathfrak g}$ and $\nu_i\in {\mathfrak g}^*$, $i=1, 2$  and we have used the previous identifications. Observe that we are identifying the elements of $T_{\alpha_g}T^*G$ with the pairs $(\xi, \nu)\in {\mathfrak g}\times {\mathfrak g}^*$.

Therefore, given the Hamiltonian $H: T^*G\equiv G\times {\mathfrak g}^*\longrightarrow \R$, we compute
\begin{equation}
dH_{(g, \alpha)}(\xi_2, \nu_2)=\langle  {\mathcal L}_g^*\left(\frac{\delta H}{\delta g}(g, \alpha)\right), \xi_2\rangle+ \langle \nu_2, \frac{\delta H}{\delta \alpha}(g, \alpha)\rangle\; , \label{hami-0}
\end{equation}
since $\frac{\delta H}{\delta \alpha}(g, \alpha)\in {\mathfrak g}^{**}={\mathfrak g}$.

We now derive the Hamilton's equations which are satisfied by the integral curves of the Hamiltonian vector field $X_H$ on $T^*G$. After left-trivialization,  $X_H(g, \alpha)=(\xi_1, \nu_1)$ where $\xi_1\in {\mathfrak g}$ and $\nu_1\in {\mathfrak g}^*$ are elements to be determined using the Hamilton's equations
\[
i_{X_H}\omega_G=dH\; .
\]
 Therefore, from expressions (\ref{omega}) and (\ref{hami-0}), we deduce that
\begin{eqnarray*}
\xi_1&=&\frac{\delta H}{\delta \alpha}(g, \alpha)\; ,\\
\nu_1&=&-{\mathcal L}_g^*\left(\frac{\delta H}{\delta g}(g, \alpha)\right)+ad_{\xi_1}^*\alpha\; .
\end{eqnarray*}
In other words, taking $\dot{g}=g\xi_1$ we obtain the {\bf Euler-Arnold equations}:
\begin{eqnarray}
\dot{g}&=&T_e{\mathcal L}_g(\frac{\delta H}{\delta \alpha}(g, \alpha))\equiv g \frac{\delta H}{\delta \alpha}(g, \alpha)\; ,\\
\dot{\alpha}&=&-{\mathcal L}_g^*\left(\frac{\delta H}{\delta g}(g, \alpha)\right)+ad_{\frac{\delta H}{\delta \alpha}(g, \alpha)}^*\alpha\; . \label{elth}
 \end{eqnarray}
If the Hamiltonian  is left-invariant, that, is there exists $h: {\mathfrak g}^*\to \R$ such that $H(g, \alpha)=H(e, \alpha)=h(\alpha)$ then we deduce that:
\begin{eqnarray*}
\dot{g}&=& g \frac{\delta h}{\delta \alpha}\, ,\\
\dot{\alpha}&=&ad_{{\delta h}/{\delta \alpha}}^*\alpha\; .
\end{eqnarray*}
The last equation is known as the {\bf Lie-Poisson equation} for a
Hamiltonian $h: {\mathfrak g}^*\to \R$.

Of course, both  Equation (\ref{elt}) and Equation (\ref{elth}) are related by the Legendre transformation when the Lagrangian $L$ is regular. This Legendre transformation is given by
\[
\begin{array}{rrcl}
Leg_L: &G\times {\mathfrak g}&\longrightarrow& G\times {\mathfrak g}^*\\
              & (g, \xi)&\longmapsto &(g, {\delta L}/{\delta \xi})
\end{array}
\]

%It is well known that both brackets are induced by reduction of the standard Lie bracket on $T^*G$ by right or left-%translation.
As we know ${\mathfrak g}^*$  is equipped with the Lie-Poisson bracket $\{\; ,\; ,\}$ which   exactly corresponds to the reduced bracket by standard  Poisson reduction from
\[
\pi: (T^*G, \omega_G)\longrightarrow (T^*G/G\equiv {\mathfrak g}^*, \{\; , \; \})
\]
where $\pi(\mu_g)=[\mu_g]\equiv T_e^*{\mathcal L}_g(\mu_g)$.

Given $\alpha\in {\mathfrak g}^*$ then it is defined the codjoint orbit by
\[
{\mathcal O}_{\alpha}:=\left\{ Ad^*_{g^{-1}}\alpha\; |\; g\in G\right\}\subseteq {\mathfrak g}^*
\]
If $t\rightarrow \alpha(t)$ is the solution of  the initial value problem
$\dot{\alpha}=ad_{{\delta h}/{\delta \alpha}}^*\alpha$ with $\alpha(0)=\alpha_0$ then we can define the curve $t\rightarrow g(t)\in G$ as the unique solution of the following first order system of differential equations: 
\[
\frac{dg}{dt}=g(t)\frac{\delta h}{\delta \alpha}(\alpha(t)), \quad g(0)=e
\]
Then
\[
\frac{d}{dt}\left( Ad^*_{g(t)^{-1}}\alpha(t)  \right)=Ad^*_{g(t)^{-1}}\left( \frac{d\alpha}{dt}(t)-ad^*_{\frac{\delta h}{\delta \alpha}(\alpha(t))}\alpha(t)\right)=0
\]
and we deduce that 
\[
 \alpha(t)\in {\mathcal O}_{\alpha(0)}
\]
given a preservation property of the continuous system that it would be important to consider in a numerical method. Therefore, since we know that coadjoint orbits are conserved quantities of the Lie-Poisson system it is natural to look for geometric integrators of the form
\[
\mu_{k+1} = Ad^*_{g_k}\mu_k\; .
\]
where $g_k$ is an appropriate element  of the Lie group $G$ (see \cite{MR1860719}). 

It $\Pi$ is the Lie-Poisson structure on ${\mathfrak g}^*$ , it is easy to check that 
\[
\begin{array}{rcl}
\sharp^{\Pi}: T^\ast{\mathfrak g}^\ast\equiv {\mathfrak g}^\ast\times {\mathfrak g}&\longrightarrow& T{\mathfrak g}^\ast \equiv {\mathfrak g}^\ast\times {\mathfrak g}^\ast\\
(\mu, \xi)&\longmapsto& (\mu, ad^*_{\xi}\mu)
\end{array}
\]

Given a hamiltonian function $h: {\mathfrak g}^*\rightarrow {\mathbb R}$ we derive the equations of motion by
the equations
\begin{equation}\label{plo}
\dot{\mu}(t)=\sharp^{\Pi}(dh(\mu(t)))
\end{equation}
or, in other words, 
\[
\dot{\mu}(t)= ad^*_{\partial h/\partial \mu}\mu(t)\, .
\]
It is well known that the flow $\Psi_t: {\mathfrak g}^*\rightarrow {\mathfrak g}^*$  of $X_h$ verifies some properties: 
\begin{enumerate}
\item It preserves the linear Poisson bracket, that is
\[
\{f\circ \Psi_t, g\circ \Psi_t\}=\{f, g\}\circ \Psi_t, ,\qquad f, g\in C^{\infty}({\mathfrak g}^*)\,.
\]
\item It preserves the hamiltonian
\[
h\circ \Psi_t=h\,.
\]
\item If all the coadjoint orbits are connected, Casimir functions are also preserved along each coadjoint orbit. 
\end{enumerate}

\section{Lie-Poisson geometric integrators} 

In this section, we will present some of the available techniques to find numerical integrators preserving as much as possible the relevant geometric structure associated to the Lie-Poisson bracket (symplectic foliation, Lie-Poisson bracket, hamiltonian, etc). 

\subsection{Discrete Lagrangian formalism}

(See for more details \cite{MR1120138,MR967725,MR1726670,MR2663389,MMM06Grupoides}).

Fixed $g\in G$, we define the set of admissible pairs
\[
C_g^2=\{(g_1, g_2)\in G\times G\; |\; g_1g_2=g\}\; .
\]
A tangent vector to the manifold $C_g^2$ is a tangent vector at $t=0$ of a curve in $C_g^2$
\[
t\in (-\epsilon, \epsilon)\subseteq \R\longrightarrow (c_1(t), c_2(t))
\]
where $c_i(t)\in G$, $c_1(t)c_2(t)=g$ and $c_1(0)=g_1$ and $c_2(0)=g_2$
All this type of curves are given by
\begin{equation}\label{curves}
c(t)=(g_1h(t), h^{-1}(t)g_2)
\end{equation}
for an arbitrary $h(t)\in G$ with $t\in (-\epsilon, \epsilon)$ and $h(0)=e$, where $e$ is the neutral element of $G$.

Fixed a discrete lagrangian $l_d: G\rightarrow \R$, we define the {\bf discrete action sum} by
\[
\begin{array}{rrcl}
S_{l_d}:&C^2_g&\longrightarrow&\R\\
               &(g_1, g_2)&\longmapsto& l_d(g_1)+l_d(g_2)
\end{array}
\]

\begin{definition}{\bf Discrete hamilton's principle}
Given $g\in G$, then $(g_1, g_2)\in C_g^2$ is a solution of the discrete lagrangian system determined by $l_d: G\rightarrow \R$ if and only if $(g_1, g_2)$ is a critical point of $S_{l_d}$.
\end{definition}

We characterize the critical points  using the curves defined in (\ref{curves}) as follows
\begin{eqnarray*}
0&=& \frac{d}{dt}\Big|_{t=0}S_{l_d}(c(t))\\
  &=& \frac{d}{dt}\Big|_{t=0}\left( l_d(g_1h(t))+l_d(h(t)^{-1}g_2)\right)\\
  &=& d(l_d\circ {\mathcal L}_{g_1})(e)(\xi)- d(l_d\circ {\mathcal R}_{g_2})(e)(\xi)
\end{eqnarray*}
where $\xi=\dot{h}(0)$. 

Alternatively, we can write these equations as follows
\begin{equation}\label{dep}
0=\lvec{\xi}(g_1)(l_d)-\rvec{\xi}(g_2)(l_d)\; , \quad \forall \xi \in {\mathfrak g}
\end{equation}
which are called {\bf discrete Euler-Poincar\'e equations}. Here $\lvec{\xi}(g)=T_e{\mathcal L}_g(\xi)$ and   $\rvec{\xi}(g)=T_e{\mathcal R}_g(\xi)$ are the left- and right-invariant vector fields. 

Also it is possible to define the two Legendre transformations by 
$Fl_d^-:  G\rightarrow {\mathfrak g}^*$ and $Fl_d^+:  G\rightarrow {\mathfrak g}^*$ by
\begin{eqnarray*}
Fl_d^-(g)&=&{\mathcal L}_g^* dl_d(g)\\
Fl_d^+(g)&=&{\mathcal R}_g^* dl_d(g)
\end{eqnarray*}
So, if we define
\[
\mu_k=Fl_d^+(g_k)={\mathcal R}_{g_k}^* dl_d(g_k)
\]
 then Equation (\ref{dep}) are equivalent to
\[
\mu_{k+1}=Fl_d^+(g_{k+1})=Fl_d^-(g_k)=Ad^*_{g_k}\mu_k
\]
which in this case are called {\bf discrete Lie-Poisson equations}. Then, it is defined an implicit map $\mu_k\mapsto \mu_{k+1}$ preserving the Lie-Poisson structure. If the discrete Lagrangian function $l_d: G\rightarrow {\mathbb R}$ is regular, that is, the Legendre transformation $Fl_d^-:  G\rightarrow {\mathfrak g}^*$ is a local diffeomorphism (or, equivalently, $Fl_d^+:  G\rightarrow {\mathfrak g}^*$ is a local diffeomorphism) then the implicit map $\mu_k\mapsto \mu_{k+1}$ is, in fact, an explicit map.

To obtain a numerical integrator for the dynamics determined by a continuous lagrangian $l: {\mathfrak g}\rightarrow {\mathbb  R}$ it is necessary to know  how closely the trajectory of the proposed  numerical method  matches the exact trajectory of the Euler-Poincar\'e equations. For variational integrators, an important tool for simplifying the error analysis is to alternatively study how closely a discrete Lagrangian matches the exact discrete Lagrangian defined by $l: {\mathfrak g}\rightarrow {\mathbb  R}$. In our case, the exact lagrangian is given by 
\[
{\mathfrak l}_h^{e}(g) = \int_{0}^{h} {\mathfrak l}(\xi_g(t))dt, \; \; \; \mbox{ for } g \in U,
\]
where $\xi_g: I \subseteq \mathbb{R} \to {\mathfrak g}$ is the unique solution of the Euler-Poincar\'e equations for ${\mathfrak l}: {\mathfrak g} \to \mathbb{R}$ such that the corresponding solution $(g, \dot{g}): I \subseteq \mathbb{R} \to TG$ of the Euler-Lagrange equations for $L(g, \dot{g})=l(g^{-1}g)$ satisfies
\[
g(0) = {\mathfrak e}, \; \; \; g(h) = g.
\]    
In \cite{MMM3} we show that if the we take as a discrete Lagrangian an approximation of order $r$ of the exact discrete Lagrangian then the associated  discrete evolution operator is also of order $r$, that is, the derived  discrete squeme is an approximation of  the continuous flow of order $r$.

\begin{example}{\bf Discrete Rigid Body Equations.}
{\rm 
Given $\Omega\in {\mathfrak so}(3)$ consider the continuous Lagrangian
\[
l(\Omega)=\frac{1}{2}\Tr(\Omega J\Omega^T)
\]
where $J$ is  positive definite matrix.
The reconstruction equation is $\dot{R}=R\Omega=R^T\Omega$ where $R	\in SO(3)$.

A discretization of this Lagrangian is given by  the discrete lagrangian $l_d: SO(3)\rightarrow {\mathbb R}$: 
\[
l_d(g_k)=l\left( \frac{g_k-I}{h}\right)=\frac{1}{2}\Tr\left( \left(\frac{g_k-I}{h}\right)J\left(\frac{g_k-I}{h}\right)^T\right)
\]
where we take the following approximation
\begin{eqnarray*}
\Omega_k&=&R_k^{-1}\dot{R}_k=R_k^T\dot{R}_k\approx R_k^T\left( \frac{R_{k+1}-R_k}{h}\right)\\
&=& \frac{R_k^TR_{k+1}-I}{h}= \frac{g_k-I}{h}
\end{eqnarray*}
and $g_k=R_k^TR_{k+1}$.

To find the critical  points of $S_{l_d}$ is equivalent to extremize $S_{\tilde{l}_d}$ where 
$\tilde{l}_d(g_k)=\Tr(g_kJ)$.
Applying discrete variational calculus, we have that
\begin{eqnarray*}
0&=&\frac{d}{dt}\Big|_{t=0}\Tr (g_kh(t)J)+\Tr (h(t)^{-1}g_{k+1}J)\\
&=&\Tr (\dot{h}(0)(Jg_k+Jg_{k+1}^T))
\end{eqnarray*}
for all $h: I\rightarrow SO(3)$ with $h(0)=I$. 

Therefore, then discrete Euler-Poincar\'e equations are
\[
Jg_k+Jg_{k+1}^T=g_k^TJ+g_{k+1}J
\]

Also it is possible to write the rigid body equations as discrete Lie-Poisson equations (see \cite{MR1120138}):
\[
\mu_{k+1}=Ad_{g_k}^*\mu_k=g_k^T \mu_k g_k
\]
where $\mu_k=Jg_k^T-g_k J$.  
}
\end{example}

\subsection{Methods based on Generating functions}

In this section we introduce the use of generating functions based on approximations of the  Hamilton-Jacobi equation to produce Lie-Poisson integrators. 
Let $G$ be a Lie group and $\mathfrak{g}$ its Lie algebra.  We know than the cotangent bundle $\tau_G: T^*G\rightarrow G$ is equipped with two natural maps $J_R: T^*G\rightarrow\mathfrak{g}^*$ and $J_L: T^*G\rightarrow\mathfrak{g}^*$. 

Let $H: {\mathfrak g}^*\rightarrow {\mathbb R}$ be a hamiltonian and 
 let $S:\mathbb{R}\times G\rightarrow
\mathbb{R}$ be a function ( a {\bf generating function})  such that the following conditions hold (see \cite{MR2839393,MR967725,MR3686003}):
\begin{enumerate}
\item {\bf Hamilton--Jacobi equation:} $\displaystyle\frac{\partial
    S}{\partial t}(t,g)+H(J_L\circ dS_t)=0$, where $S_t$ is defined
  by $S_t(g)=S(t,g)$.
\item {\it Non-degeneracy condition:} let $\xi_a$ be a basis of
  $\mathfrak{g}$. Then we assume that
  $\lvec{\xi}_a(\rvec{\xi}_b(S_t))$ is a regular matrix.
  \end{enumerate}

The idea of geometric integrators  based on generating functions  for a hamiltonian system determined by $H: {\mathfrak g}^*\rightarrow {\mathbb R}$ is to consider approximations to the generating function.
First, given the identity  element $e\in G$ take  the set $\pi_G^{-1}(e)$ which is called the set of identities of $T^*G$. It is easy to show that $\pi_G^{-1}(e)$ is a Lagrangian submanifold of $(T^*G, \omega_G)$.   Take local coordinates in $G$ around the identity $\pi_G^{-1}(e)$, say $(g^i)$, $i=1,\ldots,n$ and let $(g^i,p_i)$ be the associated natural coordinates on $T^*G$. Assume  that $g^i(e)=0$. In those coordinates 
\[
\pi_G^{-1}(e)=\{(0,p_i)\textrm{ such that }p_i\in\mathbb{R}\}.
\]

 We express the Hamilton--Jacobi equation in those coordinates: 
\[
\displaystyle\frac{\partial S}{\partial t}+H\left(J_L\circ
\frac{\partial S}{\partial g}\right)=0
\]
Now, we  approximate the solution taking the Taylor series in $t$ of $S$
  up to order $k$, $S(t,g)=\sum\limits_{i=0}^k S_i(t,g)t^i/i!+\mathcal{O}(t^{k+1})$, where
  $S_0$ is the generating function of the identity (see \cite{MR1054575,MR3686003}).
Now,  the equations for the $S_i$, $i\geq 1$ can be solved
  recursively. For instance,  we get for
  the three first terms
\begin{itemize}
\item ${S_0}(t, p_i)=0$.
\item ${S_1}(t, p_i)+H(\displaystyle\frac{\partial S_0}{\partial
    p_i},p_i)=0$.
\item ${S_2}(t, p_i)+\displaystyle\frac{\partial H}{\partial
    t}(\frac{\partial S_0}{\partial p_i},p_i)+\frac{\partial
    H}{\partial g^i}(\frac{\partial S_0}{\partial
    p_j},p_i)\frac{\partial S_1}{\partial p_i}=0$.
\end{itemize}
Each term can be obtained from the previous one by differentiating with respect
to $t$ and evaluating at $t=0$. 
Taking  all the terms obtained up to order $k$,
  $S^k=\sum\limits_{i=0}^kS_{i}t^i/i!$, we get an approximation of the
  solution of the Hamilton--Jacobi equation. It is easy to see that the
  transformation induced implicitly $\mu_1\rightarrow
  \mu_2$, $\mu^i\in {\mathfrak g}^*$ by 
\[
\begin{array}{cc}
J_R\circ\displaystyle dS^k(t,g)=\mu_1\, , & J_L\circ dS^k(t,g)=\mu_2
\end{array}
\] 
that  gives that the transformation is an approximation of order $k$ of the flow and so the
numerical method  is of order
$k$.

\subsection{Collective integrators}

The geometric integrators previously proposed are general but in some cases are  extremely complicated since they usually involve solving implicit equations in Lie groups, or using an excessive number of degrees of freedom or to calculate a computationally expensive number of 
 derivatives. In  \cite{MR3215845,MR3335215}, the authors propose the use of collective integrators. 

Let $(P, \Pi)$ be a Poisson manifold and $H: P\rightarrow {\mathbb R}$ a Hamiltonian function. 
 Let $J: M\rightarrow P$  be a  a realization of $P$ where $(M, \omega)$ is a symplectic manifold.  The function  $H\circ J: M\rightarrow {\mathbb R}$   is called the  {\bf collective Hamiltonian}. We say that a map  $\Psi: M\rightarrow M$  is {\bf collective} if there is a map $\varphi: P\rightarrow P$ such that $J\circ \Psi=\varphi\circ J$. 
 
 A {\bf collective symplectic integrator} for a hamiltonian system $H: P\rightarrow {\mathbb R}$ on a Poisson manifold $(P, \Pi)$ consist of a  full realization of $P$,  by a symplectic manifold $(M, \omega)$, $J: M\rightarrow P$  together with a symplectic integrator for $H\circ J$  that descends to a Poisson  integrator for $H$.

 For instance in \cite{MR3335215} the authors consider the map $ J: T^*{\mathbb R}^2\rightarrow {\mathbb R}^3$ defined as
\[
(q_1,q_2, p_1,p_2)\rightarrow \frac{1}{4}(2q_1q_2+2p_1p_2, 2q_1p_2-2q_2p_1, q_1^2+p^2_1-q^2_2-p^2_2),
\]
which can be used to lift a Hamiltonian function $H: {\mathbb R}^3\rightarrow {\mathbb R}$  to a collective Hamiltonian $H\circ J:T^*{\mathbb R}^2\rightarrow {\mathbb R}$ with its  canonical symplectic structure. 
To obtain a Lie-Poisson integrators, it is only necessary to  integrate the Hamiltonian vector field $X_{H\circ J}$  with a symplectic Runge-Kutta method and  to use the map $J$ to project the result back  to obtain a integrator of $X_H$ which preserve the Lie-Poisson structure on ${\mathbb R}^3$ defined by the Lie-Poisson bracket on $\mathfrak{so}(3)^*\equiv {\mathbb R}^3$.

\subsection{A geometric construction of Lie-Poisson equations from a hamiltonian function}

In this section, we propose a new possibility of constructing Lie-Poisson integrators based on the continuous hamiltonian function and a retraction map. 

Let $M_1$ and $M_2$ be $n$-dimensional manifolds and $f: M_1\rightarrow M_2$ a diffeomorphism, 
then there is a natural diffeomorphism $f_\sharp: T^*M_1\rightarrow T^*M_2$, the \emph{cotangent  lift} defined by
\[
f_\sharp(x_1, \alpha_1)=(x_2, \alpha_2)\; , \hbox{ with } \left\{
\begin{array}{rcl}
x_2&=&f(x_1)\in M_2\\
\alpha_1&=&(T_{x_1}f)^*\alpha_2\in T^*_{x_1}M_1
\end{array}
\right.
\]
for all $(x_1, \alpha_1)\in T_{x_1}^*M_1$ where  $\langle(T_{x_1}f)^*\alpha_2, v_{x_1}\rangle=\langle \alpha_2, T_{x_1}f(v_{x_1})\rangle$.

\begin{proposition} The cotangent lift $f_\sharp: T^*M_1\rightarrow T^*M_2$ of a diffeomorphism 
$f: M_1\rightarrow M_2$ is a symplectomorphism for the symplectic manifolds $(T^*M_1, \omega_{M_1})$ and $(T^*M_2, \omega_{M_2})$, in other words, 
\[
f_\sharp^*\omega_{M_2}=\omega_{M_1}
\]
\end{proposition}

However, there exist more general symplectomorphisms than the cotangent lift, for instance, the translations along the cotangent lifts on $(T^*M, \omega_{M})$. 
More generally, given a smooth function $f: M\rightarrow {\mathbb R}$ then the map $\tau_f: T^*M\rightarrow T^*M$ defined by 
\[
\tau_f (x, p)=(x, p+df(x))\, 
\]
is also a symplectomorphism.

Consider now  the cartesian product of two cotangent bundles
$
T^*M_1\times T^*M_2
$
equipped with the twisted 2-form $\Omega_{12}=-\hbox{pr}_1^*\omega_{M_1}+\hbox{pr}_2^*\omega_{M_2}$, where ${\rm pr}_i$ is the projection from $T^*M_1\times T^*M_2$ onto $T^*M_i$, $i=1, 2$. 

As we have commented any Lagrangian submanifold ${\mathcal L}$ of $T^*M_1\times T^*M_2$ being the graph of a diffeomorphism 
$F: M_1\rightarrow M_2$ guarantees that $F$ is a symplectomorphism. 
This determines a method for producing symplectomorphisms between two manifolds. 

Observe that 
the mapping 
\[
\begin{array}{rrcl}
\Phi:& T^*(M_1\times M_2)&\longrightarrow &T^*M_1\times T^*M_2\\
             & (q_1, q_2; \alpha_1, \alpha_2)&\longmapsto& (q_1, -\alpha_1; q_2, \alpha_2)
             \end{array}
             \]
is a symplectomorphism between $(T^*(M_1\times M_2), \omega_{M_1\times M_2})$     and   
$    (T^*M_1\times T^*M_2 , \Omega_{12})$.  
        
 As is known in the literature, the most  typical way to obtain a Lagrangian submanifold of $T^*(M_1\times M_2)$ consists of taking the image of ${\rm d}G$  for a function $G: M_1\times M_2\rightarrow {\mathbb R}$. Then, the Lagrangian submanifold is given by
 \[
 L_G=\{ (q_1, q_2; dG_{(q_1, q_2)})\; |\; (q_1, q_2)\in M_1\times M_2\}=\hbox{Im}\; dG
 \]
 Introducing the notation
 \begin{eqnarray*}
 dG(q_1, q_2)\doteq (D_1G(q_1, q_2), D_2 G (q_1, q_2))\in T^*_{q_1} M_1\times T^*_{q_2}M_2,
 \end{eqnarray*}
where $D_i$ denotes the derivative with respect to $q_i$, 
 then it is easy to see that 
 \[
 L_G^{\Phi}\doteq \Phi(L_G)=\{(( -D_1G(q_1, q_2);  D_2G(q_1, q_2))\in T^*M_1\times T^*M_2\}
 \]
 is a Lagrangian submanifold of $ (T^*M_1\times T^*M_2 , \Omega_{12})$. 
 
Consider a  Lagrangian submanifold  ${\mathcal L}$ of $(T^*M_1\times T^*M_2, \Omega_{12})$ which is  the graph of a symplectomorphism
 $F:  T^*M_1\rightarrow T^*M_2$, that is, ${\mathcal L}=\hbox{graph} F$, then if there exists a function $G$ on $M_1\times M_2$ such that  $L_G^{\Phi}={\mathcal L}$ we will say that  $F$ is the symplectomorphism generated by $G$ and $G$ is 
 the {\bf generating function} of $F$.
 
 General Lagrangian submanifolds of $(T^*M_1\times T^*M_2, \Omega)$ are called {\bf canonical relations}
  when they are thought as generators of implicit symplectomorphisms between $T^*M_1$ and $T^*M_2$.

%%If $L_{12}$ is a Lagrangian submanifold  or canonical relation of $(T^*M_1\times T^*M_2 , \Omega_{12})$ and $L_{23}$ is a canonical relation of 
%%$(T^*M_2\times T^*M_3 , \Omega_{23})$, then we can form the \emph{composition}  
%%\[
%%L_{13}=L_{12}\circ L_{23}\in T^*M_1\times T^*M_3
%\]
%defined by
%\begin{eqnarray}
%L_{13}&=&\{(q_1, \alpha_1; q_3, \alpha_3)\in T^*M_1\times T^*M_3\; \mid\;  \exists (q_2, \alpha_2)\in T^*M_2, \hbox{  such that  } \label{tyu}\\
%&&(q_1, \alpha_1, q_2, \alpha_2)\in T^*M_1\times T^*M_2\ \hbox{and}\ (q_2, \alpha_2, q_3, \alpha_3)\in T^*M_2\times T^*M_3\}\nonumber
%\end{eqnarray}
%Denote by $\Delta_{T^*M_2}=\{(\mu_2, \mu_2)\in T^*M_2\times T^*M_2\; \mid\; \mu_2\in T^*M_2\}$
%\begin{theorem}
%If the canonical relations $L_{12}$ and $L_{23}$ 
%intersect cleanly,  that is, 
%\[
%(pr_2(L_{12}), pr_1(L_{23}))  \pitchfork\Delta_{T^*M_2}
%\]
%then their composition $L_{13}$ is an immersed Lagrangian
%submanifold of 
% $(T^*M_1\times T^*M_3 , \Omega_{13})$.
%\end{theorem}
%
%
%
%If $L_{12}=L_{G_{12}}^\sigma$ and $L_{23}=L_{G_{23}}^\sigma$ for some functions $G_{12}: M_1\times M_2\rightarrow {\mathbb R}$ and $G_{23}: M_2\times M_3\rightarrow {\mathbb R}$
%then if we define for any $(q_1, p_1)$ and $(q_3, p_3)$ the map
%\[
%q_2\longrightarrow G_{12}(q_1,q_2)+G_{23}(q_2, q_3)
%\]
%and assume that this map has a unique nondegenerate critical point $q_2^*$, then 
%the map
%\[
%G_{13}(q_1, q_3)= G_{12}(q_1,q^*_2)+G_{23}(q^*_2, q_3)
%\]
%is a generating function of the Lagrangian submanifold $L_{13}$ defined on  (\ref{tyu}). 

Let  $({\mathfrak g}^*, \Pi)$ be a linear Poisson bracket on the dual of a Lie algebra ${\mathfrak g}$, then given a Hamiltonian $H: {\mathfrak g}^*\rightarrow {\mathbb R}$ we have the  corresponding Lie-Poisson equations: 
\begin{equation*}\label{plo}
\dot{\mu}(t)=\sharp^{\Pi}(dH(\mu(t)))
\end{equation*}
These equations define the flow of the hamiltonian vector field $X_h$: 
$$\Psi_t: {\mathfrak g}^*\rightarrow {\mathfrak g}^*\; ,$$ 
which obvioulsy  is a Poisson morphism or, equivalently, 
 ${\rm Graph}\, \Psi_t$ is a coisotropic submanifold of 
 $(\overline{{\mathfrak g}^*}\times {\mathfrak g}^*,\Pi_{\overline{{\mathfrak g}^*}\times {\mathfrak g}^*})$. 

A {\bf geometric integrator} for the Hamiltonian vector field $X_h$  consists of a coisotropic submanifold $(\overline{{\mathfrak g}^*}\times {\mathfrak g}^*,\Pi_{\overline{{\mathfrak g}^*}\times {\mathfrak g}^*})$ ``near of" ${\rm Graph}\, \Psi_t$.

\begin{example}
{\rm 
Our objective is to generate a Lagrangian submanifold ${\mathcal L}$ of $(T^*G, \omega_G)$. We will obtain this Lagrangian submanifold from a hamiltonian 
$H: {\mathfrak g}^*\rightarrow {\mathbb R}$ and a retraction map 
  $\tau: {\mathfrak g}\to G$, which is an analytic local
diffeomorphism around the identity such that $\tau(\xi)\tau(-\xi)={e}$,
where $\xi\in\mathfrak g$. Thereby $\tau$ provides a local chart on the Lie
group. 

Given a  retraction map $\tau$, the right trivialized tangent $\mbox{d}\tau_{\xi}:\mathfrak{g}\rightarrow\mathfrak{g}$ and the  inverse $\mbox{d}\tau_{\xi}^{-1}:\mathfrak{g}\rightarrow\mathfrak{g}$ are defined for all   $\eta\in\mathfrak{g}$ as follows
  \begin{eqnarray*}
    T_{\xi}\tau (\eta)&=&T_{ e} {\mathcal R}_{\tau(\xi)}\left(\mbox{d}\tau_{\xi}(\eta)\right)\equiv \mbox{d}\tau_{\xi}(\eta)\,\tau(h\xi),\\
    T_{\tau(\xi)}\tau^{-1}((T_{ e} {\mathcal R}_{\tau(\xi)})\eta)&=&\mbox{d}\tau^{-1}_{\xi}(\eta).
  \end{eqnarray*}

Consider as a particular example
$\tau: {\mathfrak g}\rightarrow G$ defined by
\[
\tau(\xi)=\hbox{exp } (\frac{h\xi}{2})
\]
Then, $\tau_{\sharp}:  T^*{\mathfrak g}\equiv {\mathfrak g}\times {\mathfrak g}^*\rightarrow T^*G$ is a symplectomorphism. It is given by
\begin{eqnarray*}
\tau_{\sharp}(\xi, \mu)&=& (g=\hbox{exp } (\frac{h\xi}{2}), (T_{g} \hbox{exp}^{-1})^*\mu )\\
&=& (g=\hbox{exp } (\frac{h\xi}{2}), \frac{h}{2}{\mathcal R}_{g^{-1}}^*((d_{h\xi/2} \hbox{exp}^{-1})^*\mu) )
\end{eqnarray*}

Consider now the following full realization
\[
\begin{array}{rcl}
J:&T^*G&\rightarrow \overline{ {\mathfrak g}^*}\times {\mathfrak g}^*\\
   &(g, \eta)&\longmapsto ({\mathcal R}_g^*\eta, {\mathcal L}_g^*\eta)
   \end{array}
\]
Therefore, 
\[
J\circ\tau_{\sharp}(\xi, \mu)=( \frac{h}{2}(d_{h\xi/2} \hbox{exp}^{-1})^*\mu), \frac{h}{2}\hbox{Ad}_{g}^*((d_{h\xi/2} \hbox{exp}^{-1})^*\mu))
   \]
where $g=\hbox{exp } (\frac{h\xi}{2})$.  

\

{\bf Lagrangian mechanics}. 

Given a lagrangian $l: {\mathfrak g}\rightarrow {\mathbb R}$,
then $\hbox{Im }dl$ is a Lagrangian submanifold of $T^*{\mathfrak g}$: 
\[
\hbox{Im }dl=\{ (\xi, \frac{\delta l}{\delta \xi}(\xi))\; |\; \xi\in {\mathfrak g} \}
\]
 Therefore, 
\[
J\tau_{\sharp}(\hbox{Im }dl)
\]
is a coisotropic submanifold of  $\overline{{\mathfrak g}^*}\times {\mathfrak g}^*$. 
This coisotropic submanifold  implicitly defines the  Poisson map $\mu_k\rightarrow \mu_{k+1}$ defined by
\begin{eqnarray*}
\mu_k&=&\frac{h}{2}(d_{h\xi/2} \hbox{exp}^{-1})^*\frac{\delta l}{\delta \xi}\\
\mu_{k+1}&=&
\frac{h}{2}\hbox{Ad}_{g}^*((d_{h\xi/2} \hbox{exp}^{-1})^*\frac{\delta l}{\delta \xi})
\end{eqnarray*}
where $g=\hbox{exp } (\frac{h\xi}{2})$. Therefore, the Lie-Poisson geometric integrator is
\begin{eqnarray*}
\mu_k&=&\frac{h}{2}(d_{h\xi/2} \hbox{exp}^{-1})^*\frac{\delta l}{\delta \xi}\\
\mu_{k+1}&=&\hbox{Ad}^*_g\mu_k\\
g&=&\hbox{exp } (\frac{h\xi}{2})
\end{eqnarray*}

\

{\bf Hamiltonian mechanics}.

Given a hamiltonian  $H: {\mathfrak g}^*\rightarrow {\mathbb R}$ not necessarily regular.
Then $\hbox{Im }dH$ is a Lagrangian submanifold of $T^*{\mathfrak g}^*$. 
Using the canonical antisymplectomorphism between $T^*{\mathfrak g}^*={\mathfrak g}^*\times {\mathfrak g}$ and $T^*{\mathfrak g}={\mathfrak g}\times {\mathfrak g}^*$
given by 
$\varphi(\mu, \xi)\rightarrow (\xi, \mu)$. 
\[
\varphi(\hbox{Im }dH)=\{  (\frac{\delta H}{\delta \mu}(\mu), \mu)\; |\; \mu\in {\mathfrak g}^* \}
\]
 Therefore, 
\[
J\tau_{\sharp}(\varphi(\hbox{Im }dH))
\]
is also a  coisotropic submanifold of  $\overline{{\mathfrak g}^*}\times {\mathfrak g}^*$. 
This coisotropic submanifold  implicitly defines the  Poisson map $\mu_k\rightarrow \mu_{k+1}$ defined by
\begin{eqnarray*}
\mu_k&=&\frac{h}{2}(d_{\frac{h}{2}\frac{\delta H}{\delta \mu}} \hbox{exp}^{-1})^*\mu)\\
\mu_{k+1}&=&
\frac{h}{2}\hbox{Ad}_{g}^*((d_{\frac{h}{2}\frac{\delta H}{\delta \mu}} \hbox{exp}^{-1})^*\mu)
\end{eqnarray*}
where $g=\hbox{exp } (\frac{h}{2}\frac{\delta H}{\delta \mu})$.
 Therefore, the Lie-Poisson geometric integrator is
\begin{eqnarray*}
\mu_k&=&\frac{h}{2}(d_{\frac{h}{2}\frac{\delta H}{\delta \mu}} \hbox{exp}^{-1})^*\mu)\\
g&=&\hbox{exp } (\frac{h}{2}\frac{\delta H}{\delta \mu})\\
\mu_{k+1}&=&\hbox{Ad}^*_g\mu_k
\end{eqnarray*}

}
\end{example}
\subsection{Energy-preserving integrators based on discrete gradients}

In previous sections, we have seen that given a Poisson manifold $(M, \Pi)$ and Hamiltonian function $H: M\rightarrow \R$
the integral curves of a Hamiltonian vector field are written as
\begin{equation}\label{qer}
\frac{dx}{dt}=X_H(x(t))=\sharp^{\Pi}(dH)
\end{equation}
Working in a local neighborhood $(U, (x^i, \ldots, x^m))$ we can write Equations (\ref{qer}) as
\[
\frac{dx}{dt}=\Pi(x) \nabla H(x)
\]
where $\Pi(x)=(\Pi^{ij}(x))$ represents the skew symmetric matrix given by $\Pi^{ij}(x)=\{x^i, x^j\}(x)$

Then $\overline{\nabla}H: U\times U\longrightarrow U$ is called a {\bf discrete gradient} of $H$ if it is continuous and it satisfies
\begin{eqnarray*}
\overline{\nabla}H(x,x')\cdot (x'-x)&=&H(x')-H(x)\\
\overline{\nabla}H(x,x)&=&\nabla H(x)
\end{eqnarray*}
for all $x, x'\in U$.
There are different possibilities to construct discrete gradients for a function $H: U\rightarrow {\mathbb R}$. For instance, the 
{\bf mean value discrete gradient}
$$
\overline{\nabla}H(x,x'):=\int_0^1 \nabla H ((1-\xi)x+\xi x')d\xi
$$
or the
{\bf  midpoint discrete gradient}
$$
\overline{\nabla}_H(x,x'):=\nabla H\left( \frac{1}{2}(x'+x)\right)+\frac{H(x')-H(x)-\nabla H\left( \frac{1}{2}(x'+x)\right)\cdot (x'-x)}{|x'-x|^2}(x'-x), 
$$
where  $x'\not=x$.
See \cite{MR943488,MR1694701} for more details

We define an energy-preserving integrator by
$$
\frac{x'-x}{h}=\tilde{\Pi}(x,x',h)\overline{\nabla}H(x,x')
$$
where $\tilde{\Pi}$ is an approximation of $\Pi$, that is $\tilde{\Pi}(x,x,h)=\Pi(x)$. 
This integrator preserves the Hamiltonian  but, in general, it is not preserving the Poisson structure.

\section*{Acknowledgements}
The author has been partially supported by Ministerio de Econom\'ia, Industria y Competitividad (MINEICO, Spain) under grants MTM 2013-42870-P, MTM 2015-64166-C2-2P, MTM2016-76702-P and ``Severo Ochoa Programme for Centres of Excellence'' in R\&D (SEV-2015-0554). The author would like to thank the anonymous reviewer for the  valuable comments and suggestions to improve the quality of the paper.

 \bibliography{References}

\end{document}